\documentclass{article}

\usepackage{amsmath,amssymb}
\usepackage[mathscr]{eucal}

\usepackage{mathtools}                                                      
\mathtoolsset{showonlyrefs=true}                                    

\usepackage{pdfsync}

\usepackage{dsfont}

\usepackage{eufrak}                                             

\let\counterwithin\relax

\makeatletter
\newif\if@check@engine  \@check@enginetrue
\makeatother

\usepackage{amsmath}                                                        
\usepackage{graphicx}                                                       
\usepackage{subfigure}
\usepackage{enumitem}                                                    
\usepackage{hyperref}
\usepackage{amssymb}                                                        
\usepackage[mathscr]{eucal}                                             
\usepackage{cancel}                                                             
\usepackage[normalem]{ulem}                                                                 
\usepackage{pstricks}
\usepackage{rotating}
\usepackage{lscape}
\usepackage[paperwidth=8.5in,paperheight=11in,top=1.25in, bottom=1.25in, left=1.00in, right=1.00in]{geometry}
\usepackage{mathtools}                                                      
\mathtoolsset{showonlyrefs=true}                                    
\usepackage{fixltx2e,amsmath}                                           
\MakeRobust{\eqref}
\usepackage{mathdots}
\usepackage{amsthm}                                                             
\allowdisplaybreaks                                                             
\usepackage{chngcntr}

\usepackage{ulem}       

\usepackage{titlesec}
\usepackage[titletoc,toc,title]{appendix}

\usepackage{apptools, etoolbox}

\AtAppendix{\counterwithin{theorem}{section}}

\usepackage{color}

\newtheorem{theorem}{Theorem}[section]
\newtheorem{proposition}[theorem]{Proposition}
\newtheorem{corollary}[theorem]{Corollary}
\newtheorem{lemma}[theorem]{Lemma}
\newtheorem{definition}[theorem]{Definition}

\newtheorem{remark}[theorem]{Remark}
\newtheorem{assumption}[theorem]{Assumption}

\def \aa {\mathbf{a}}
\def \bb {\mathbf{b}}
\def \cc {\mathbf{c}}
\def \ff {\mathbf{f}}
\def \gg {\mathbf{g}}
\def \hh {\mathbf{h}}

\def \LL {\mathbf{L}}
\def \Lc {\mathbf{G}}
\def \AA {\mathbf{A}}
\def \MM {\l}
\def \MMd {\l}
\def \DD {\nabla}
\def \Hh {\mathcal{H}}

\def \o {{\omega}}
\def \a {{\alpha}}
\def \b {{\beta}}
\def \d {{\delta}}
\def \l {{\lambda}}

\def \G {{\Gamma}}
\def \s {{\sigma}}
\def \w {{\omega}}

\def \R {{\mathbb {R}}}
\def \N {{\mathbb {N}}}

\def \x {{\xi}}
\def \e {{\varepsilon}}

\def \r {{\varrho}}

\def \t {{\tau}}
\def \t {{\tau}}

\def \m {{\mu}}
\def \y {{\eta}}

\def \tt {{\t}}

\def \O {{\Omega}}
\def \phi {{\varphi}}

\def \tilde {\widetilde}

\def\p{\partial}

\def \B {{\cal{B}}}

\def \o {{\omega}}
\def \a {{\alpha}}
\def \b {{\beta}}
\def \d {{\delta}}
\def \l {{\lambda}}

\def \G {{\Gamma}}
\def \s {{\sigma}}
\def \w {{\omega}}

\def \R  {{\mathbb {R}}}
\def \x {{\xi}}

\def \e {{\varepsilon}}

\def \t {{\tau}}

\def \m {{\mu}}
\def \y {{\eta}}

\def \p {{\partial}}
\def \a {{\alpha}}
\def \O {{\Omega}}

\def \d {{\delta}}

\def \o {{\omega}}
\def \a {{\alpha}}
\def \b {{\beta}}
\def \d {{\delta}}

\def \G {\Ga}

\def \Ga {{\Gamma}}

\def \Ito {It\^o }
\def \s {{\sigma}}

\def \w {{\omega}}
\def \R {{\mathbb {R}}}
\def \N {{\mathbb {N}}}

\def \x {{\xi}}
\def \e {{\varepsilon}}

\def \r {{\varrho}}

\def \t {{\tau}}
\def \t {{\tau}}

\def \m {{\mu}}
\def \y {{\eta}}

\def \tt {{\t}}

\def \O {{\Omega}}
\def \phi {{\varphi}}

\def \tilde {\widetilde}

\def\l {\lambda}

\def \F {\mathcal{F}}
\def \B {\mathscr{B}}

\def \Ã  {{\`a }}
\def \è {{\`e }}
\def \ò {{\`o }}
\def \ù {{\`u }}

\def \caratt {{\mathds{1}}}

\makeatletter
\def\section{\@startsection {section}{1}{\z@}{3.25ex plus 1ex minus
 .2ex}{1.5ex plus .2ex}{\large\bf}}
\def\subsection{\@startsection{subsection}{2}{\z@}{3.25ex plus 1ex minus
 .2ex}{1.5ex plus .2ex}{\normalsize\bf}}
\@addtoreset{equation}{section} 

\makeatother

\title{\empty}

\author{\empty}

\date{\empty}

\numberwithin{equation}{section}

\begin{document}

\title{The parametrix method for parabolic SPDEs}
\author{
Andrea Pascucci\thanks{(Corresponding author) Dipartimento di Matematica, Universit\`a di Bologna,
Bologna, Italy. \textbf{e-mail}: andrea.pascucci@unibo.it} \and Antonello
Pesce\thanks{Dipartimento di Matematica, Universit\`a di Bologna, Bologna, Italy. \textbf{e-mail}:
antonello.pesce2@unibo.it} }

\date{This version: \today}

\maketitle

\begin{abstract}
We consider the Cauchy problem for a linear stochastic partial differential equation. By extending
the parametrix method for PDEs whose coefficients are only measurable with respect to the time
variable, we prove existence, regularity in H\"older classes and estimates from above and below of
the fundamental solution. This result is applied to SPDEs by means of the It\^o-Wentzell formula,
through a random change of variables which transforms the SPDE into a PDE with random
coefficients.
\end{abstract}


\section{Introduction}
Stochastic partial differential equations (SPDEs) arise in many applications in probability theory
and
in particular in the study of filtering problems (see e.g. \cite{MR0336801}, \cite{MR3334279}). 
Assume that $(X_{t},Y_{t})$ is a diffusion  where $X$ represents a signal that is not directly
observable and we want to extract information about $X$ from $\F^{Y}_{t}=\s(Y_{s},\, s\le t)$ that
is the filtration of the observations on $Y$. Then, under natural assumptions, for any bounded and
measurable function $f$ we have
  $$E\left[f(X_{t})\mid\F^{Y}_{t}\right]=\int_{\R^{d}} f(x)p_{t}(x)dx$$
where $p_{t}(x)$ denotes the conditional density of $X_{t}$ given $\F_{t}^{Y}$: it turns out that
$p_{t}$ satisfies a SPDE of the form
\begin{equation}\label{e16}
  dp_{t}(x)=\LL_{t}p_{t}(x)dt+\Lc_{t}p_{t}(x)dW_{t},
\end{equation}
where $\LL_{t}$ is a second-order elliptic operator and $\Lc_{t}$ is a first-order operator. The
coefficients of $\LL_{t}$ and $\Lc_{t}$ may depend on $t,x,Y_{t}$ and are therefore random and
typically not smooth. A very particular case is when $Y\equiv 0$: then $\Lc_{t}= 0$ and
\eqref{e16} reduces to the classical {forward} Kolmogorov equation for the transition density
$p_{t}$ of $X_{t}$. In the general case, $p_{t}$ can be referred to as the  stochastic
fundamental solution of \eqref{e16}. 
The aim of this paper is to prove existence, regularity and estimates from above and below of
$p_{t}$ by using a classical tool from PDEs' theory, the parametrix method.

Let $(\O,\F,P)$ be a complete probability space with an increasing filtration
$\left(\F_{t}\right)_{t\ge 0}$ of complete with respect to $(\F,P)$ $\s$-fields
$\F_{t}\subseteq\F$. Let $d_1\in\N$ and let $W^k$, $k=1,\cdots, d_1$, be one-dimensional
independent Wiener processes with respect to $\left(\F_{t}\right)_{t\ge 0}$. We consider the
parabolic SPDE
\begin{equation}\label{SPDE}
  du_t(x)=\left(\LL_tu_t(x)+\ff_t(x)\right)dt+
  \Lc_{\s_t^k}u_t(x)
  dW^k_t,
\end{equation}
where $\LL_t$ is the second-order operator
\begin{equation}
 \LL_tu_t(x)=\frac{1}{2}
 \aa_t^{ij}(x)\p_{ij}u_t(x)+
 \bb_t^j(x)\p_ju_t(x)+\cc_t(x)u_t(x)
\end{equation}
and $\Lc_{\s_t^k}$ is the first-order operator
\begin{equation}
 \Lc_{\s_t^k}u_t(x)=
 \s^{ik}_t(x)\p_iu_t(x).
\end{equation}
Throughout the paper, the summation convention over repeated indices is enforced regardless of
whether they stand at the same level or at different ones. The point of $\R^{d}$ is denoted by
$x=(x_{1},\dots,x_{d})$ and we set $\p_{i}=\p_{x_{i}}$, $\DD=\left(\p_{1},\dots,\p_{d}\right)$ and
$\p_{ij}=\p_{i}\p_{j}$. A function $u=u_{t}(x,\w)$ on $[0,\infty)\times\R^{d}\times\O$ is denoted
by $u_{t}(x)$ and we shall systematically omit the explicit dependence on $\w\in\O$. The
coefficients $\aa_t$, $\bb_t$, $\cc_t$, $\ff_t$ and $\s^k_t$ 
in \eqref{SPDE} are intended to be random and not smooth. 

The Cauchy problem for evolution SPDEs has been studied by several authors. Under coercivity
conditions analogous to uniform ellipticity for PDEs, there exists a complete theory in Sobolev
spaces (see \cite{MR0651582}, \cite{MR0501350}, \cite{MR1135324}, \cite{MR1042066} and references
therein) and in the spaces of Bessel potentials (see \cite{MR1377477} and \cite{MR1661766}). 
Classical solutions in H\"older classes were first considered in \cite{MR0370764}, \cite{Shimizu}
and more recent results were proved in \cite{Chow94} and \cite{Mikulevicius}, though the authors
only considered equations with non-random coefficients and with no derivatives of the unknown
function in the stochastic term.

In the last decades, the use of analytical or PDE techniques in the study of SPDEs has become
widespread. For instance, the results in \cite{Chow94}, \cite{Mikulevicius}, \cite{Zatezalo} are
based on classical methods of deterministic PDEs, such as Duhamel principle and a priori Schauder
estimates; the $L^{p}$ estimates in \cite{Matoussi} are proved by adapting the classical Moser's
iterative argument; \cite{Sowers_mono} provides short-time asymptotics of random heat kernels. A
further remarkable example is given by the recent series of papers by Krylov \cite{Krylov14,
Krylov17, MR3334279} where the H\"ormander's theorem for SPDEs is proved; see also the very recent
results in \cite{MR3706782} for backward SPDEs.

In this paper we extend another classical tool that, to the best of our knowledge, has not yet
been considered in the study of SPDEs, the well-known {\it parametrix method} for the construction
of the fundamental solution of PDEs with H\"older continuous coefficients. There are two main
problems that one faces when trying to apply the parametrix method to SPDEs: the lack of a Duhamel
principle for SPDEs and the roughness of the coefficients, that are assumed to be only measurable
in time. In the first part of the paper, we use the It\^o-Wentzell formula \cite{Ventcel} to make
a random change of variables which transforms the SPDE in a PDE with random coefficients; the
latter admits a Duhamel principle and, in the second part of the paper, we use it to extend the
parametrix method to parabolic PDEs with coefficients measurable in the time variable.

{This paper does not pretend to encompass the most general assumptions but rather investigate the
possible use of the parametrix method in the stochastic framework; as such, it has to be intended
as a first step of a research programme aiming at considering more general classes of possibly
degenerate SPDEs. More precisely, it is very likely that the techniques used in this paper can be
applied to SPDEs satisfying the strong H\"ormander condition, such as those considered in
\cite{Krylov14}. A more challenging problem is to consider the Langevin SPDE
  $$du_{t}(x,y)=\left(\frac{a_{t}(x,y)}{2}\p_{xx}u_{t}(x,y)+x\p_{y}u_{t}(x,y)\right)dt+\s_{t}(x,y)\p_{x}u_{t}(x,y)dW_{t},\qquad (x,y)\in\R^{2}.$$
This equation has unbounded drift coefficient and satisfies the {\it weak} H\"ormander condition.
The pa\-ra\-me\-trix method has been generalized to {\it deterministic} (i.e. with $\s_{t}\equiv 0
$) Langevin PDEs in \cite{Polidoro2}, \cite{DiFrancescoPascucci2} and \cite{MR2744877}; however,
contrary to the uniformly parabolic case considered in the present paper, the intrinsic H\"older
regularity in the spatial variables cannot be studied independently from the time variable as it
was recently shown in \cite{MR3660883}. This is an additional issue that needs careful
investigation and is the subject of the forthcoming paper \cite{PascucciPesce2}.}

The paper is organized as follows. In Section \ref{sec2} we introduce the basic notations and
state our main result, Theorem \ref{t2b}; for illustrative purposes, the particular case of the
stochastic heat equation is discussed in Subsection \ref{sec21}. In Section \ref{ito} we recall
the It\^o-Wentzell formula and provide some estimate for the related flow of diffeomorphisms. In
Section \ref{parametrixsec} we present the parametrix method. Since the complete proofs are rather
technical and to a large extent similar to the classical case, we only provide the details on
those aspects that require significant modifications: in particular, in Section \ref{lowr} we
present a proof of the Gaussian lower bound for the fundamental solution which requires some non
trivial adaptation of an
original argument by Aronson (cf. \cite{FaSt}). 

\section{Assumptions and main results}\label{sec2}
Before stating our main theorems, we need to introduce some basic definitions and notations to be
used throughout the paper. Let $k\in\N$, $\a\in(0,1)$ and $0\le t<T$. Denote:
{\begin{itemize}
 \item $C^{0}_{t,T}=C^{0}_{t,T}(\R^{d})$ the space of all measurable
functions $f=f_{s}(x)$ on $[t,T]\times \R^d$ that are continuous in $x$;
  \item $C^{\a}_{t,T}=C^{\a}_{t,T}(\R^{d})$ the space of functions $f\in C^{0}_{t,T}$
  that are $\a$-H\"older continuous in $x$ uniformly with respect to $s$, that is
  $$\sup_{s\in[t,T]\atop x\neq y}\frac{|f_{s}(x)-f_{s}(y)|}{|x-y|^{\a}}<\infty;$$
  \item $C^{\a}_{t,T,\text{\rm loc}}$ the space of functions $f\in C^{0}_{t,T}$ such that
  $$\sup_{s\in[t,T]\atop x,y\in K,\, x\neq y}\frac{|f_{s}(x)-f_{s}(y)|}{|x-y|^{\a}}<\infty$$
  for every compact subset $K$ of $\R^{d}$;
   \item $C^{k}_{t,T}$  the space of all measurable
   functions $f=f_{s}(x)$ on $[t,T]\times \R^d$
    that are $k$-times differentiable w.r.t. $x$ with $\p^{\b}f\in C^{0}_{t,T}$ for any multi-index
$\b$ of height $|\b|= k$;
  \item $C^{k+\a}_{t,T}$ (resp. $C^{k+\a}_{t,T,\text{\rm
loc}}$) the space of functions $f\in C^{k}_{t,T}$ 
with $\p^{\b}f\in C^{\a}_{t,T}$  (resp. $\p^{\b}f\in C^{\a}_{t,T,\text{\rm loc}}$) for any
multi-index $\b$ of height $|\b|\le k$.
\end{itemize}}

We use boldface to denote the stochastic versions of the previous functional spaces.
More precisely, let {now} $k\in\N\cup\{0\}$, $\a\in[0,1)$ and $\mathcal{P}_{t,T}$ be the predictable
$\s$-algebra on $[t,T]\times\O$. We denote by $\mathbf{C}^{k+\a}_{t,T}$ the family of functions
$f=f_{s}(x,\o)$ on $[t,T]\times\R^{d}\times\O$ such that:
\begin{itemize}
 \item[i)] $(s,x)\mapsto f_{s}(x,\o)\in C^{k+\a}_{t,T}$ 
 {$P$-a.s.;}
 \item[ii)] $(s,\w)\mapsto f_{s}(x,\o)$ is $\mathcal{P}_{t,T}$-measurable for any $x\in\R^{d}$.
\end{itemize}
Moreover, $\mathbf{bC}^{k+\a}_{t,T}$ is the space of functions $f\in\mathbf{C}^{k+\a}_{t,T}$ such
that
  $$\sum_{|\b|\le k}\sup_{s\in[t,T]\atop x\in\R^d}|\p^{\b}f_{s}(x)|
  <\infty\qquad P\text{-a.s.}$$
We say that $f=f_{s}(x)$ is {\it non-rapidly increasing uniformly on $(t,T]\times \R^{d}$} if, for
any $\d>0$, $e^{-\d|x|^{2}}\left|f_{s}(x)\right|$ is a bounded function on $(t,T]\times \R^{d}$,
$P$-a.s.; in case $f$ does not depend on $s$, we simply say that $f$ is non-rapidly increasing on
$\R^{d}$.

\begin{definition}\label{d3b} A stochastic fundamental solution $\mathbf{\G}=\mathbf{\G}(t,x;\t,\x)$ for
the SPDE \eqref{SPDE} is a function defined for $0\le\t<t\le T$ and $x,\x\in\R^{d}$, such that for
any $(\t,\x)\in[0,T)\times\R^{d}$ we have:
\begin{itemize}
  \item[i)] {$\mathbf{\G}(\cdot,\cdot;\t,\x)\in \mathbf{C}^{2}_{t_{0},T}(\R^{d})$} and with probability
  one satisfies
\begin{equation}\label{e50}
 \mathbf{\G}(t,x;\t,\x)=\mathbf{\G}(t_{0},x;\t,\x)+
  \int_{t_{0}}^t \LL_s \mathbf{\G}(s,x;\t,\x)ds+
 \int_{t_{0}}^t\Lc_{\s_s^k}\mathbf{\G}(s,x;\t,\x)dW^k_s
\end{equation}
 for $\t<t_{0}\le t\le T$ and $x\in\R^d$;
  \item[ii)] for any continuous and non-rapidly increasing function $\phi$ on $\R^{d}$
    $$\lim_{(t,x)\to(\t,\x)\atop t>\t}\int_{\R^{d}}\mathbf{\G}(t,x;\t,y)\phi(y)dy=\phi(\x),\qquad P\text{-a.s.}$$
\end{itemize}
\end{definition}

Next we state the standing assumptions on the coefficients of the SPDE \eqref{SPDE}.
\begin{assumption}[\bf Coercivity]\label{Ass1}
Let
   $$\mathbf{A}_{t}(x):=\left(\aa^{ij}_{t}(x)-
   \sigma_{t}^{ik}(x)\sigma_{t}^{jk}(x)\right)_{i,j=1,\dots,d}.$$
There exists a positive random variable $\MM$ such that
\begin{equation}
 \langle \AA_{t}(x)\x,\x\rangle\ge \MM|\xi|^{2},\qquad
 t\in [0,T],\ x,\x\in \R^{d},\ P\text{-a.s.}
\end{equation}
\end{assumption}

\begin{assumption}[\bf Regularity]\label{Ass2}
For some $\a\in(0,1)$ and for every $i,j = 1,\dots, d$ and $k = 1,\dots,d_1$, we have
 $$\aa^{ij},\ \bb^j, \ \cc\in\mathbf{bC}^{\a}_{0,T}\quad \text{and} \quad
 \s^{ik}\in\mathbf{bC}^{3+\a}_{0,T}.$$
\end{assumption}

We now introduce a random change of coordinates that will play a central role in the following
analysis. We fix $(\t,x)\in[0,T)\times\R^{d}$ and consider the stochastic ordinary differential
equation
\begin{equation}\label{IW-SDE}
 x_{t}=x-\int_{\t}^{t}\s^{k}_{s}(x_s)dW^{k}_{s},\qquad t\in[\t,T].
\end{equation}
It is well-known (see, for instance, Theor. 4.6.5 in \cite{MR1070361})) that, under Assumption
\ref{Ass2}, equation \eqref{IW-SDE} admits a solution $X=X_{\t,t}(x,\w)$ that is a stochastic flow
of diffeomorphisms: precisely, $X_{\t,t}\in \mathbf{C}^{3+\a'}_{\t,T}$, for any $\a'<\a$, the
matrix $\DD X_{\t,t}(x)$ satisfies
\begin{align}\label{eq5}
 \DD X_{\t,t}(x)&=I_d-\int_{\t}^t \DD \s^k_s(X_{\t,s}(x))\DD X_{\t,s}(x)dW_s^k,
 \end{align}
and, for any $i,j=1,\dots,d$, 
$\p_{ij}^2 X_{\t,t}(x)$ satisfies
{\begin{equation}\label{eq5b}
 \p_{ij}^2X^h_{\t,t}(x)=-\int_{\t}^t \left((\DD \s_s^{k}(X_{\t,s}(x))\p_{ij}^2X_{\t,s}(x))_h+ \left((\DD
 X_{\t,s}(x))^{\ast}\DD^2 \s_s^{h k}(X_{\t,s}(x))\DD X_{\t,s}(x)\right)_{ij}\right)dW^k_s
\end{equation}}
with probability one.

{Since we are going to use $X$ as a global change of variables, we need some control over the
stochastic integrals in \eqref{eq5} and \eqref{eq5b} for $x$ varying in $\R^{d}$: this issue is
addressed in Section \ref{ito} (see, in particular, Proposition \ref{lemma1}) under the following
additional condition.}
\begin{assumption}\label{Ass3}
There exist $\e>0$ and a random variable $M\in L^{\bar{p}}(\O)$, with
$\bar{p}>\max\left\{2,d,\frac{d}{2\e}\right\}$, such that
 $$\sup_{s\in[0,T]\atop x\in\R^d}(1+|x|^2)^{\e}|\p^{\b}\s^k_s(x)|\leq M \qquad \text{P-}a.s.$$
for every $k=1,\cdots, d_1$ and multi-index $\b$ with $1\leq |\b|\leq 3$.
\end{assumption}
Assumption \ref{Ass3} is a rather weak condition on the first, second and third order derivatives
of $\s$: clearly, it is satisfied under the very particular cases of $\s$ constant or $\s$ with
compact support.

Our main result is the following
\begin{theorem}\label{t2b}
Let Assumptions \ref{Ass1}, \ref{Ass2} and \ref{Ass3} be in force. Then there exists a fundamental
solution $\mathbf{\G}$ for the SPDE \eqref{SPDE}. Moreover, there exist two positive random
variables $\m_{1}$ and $\m_{2}$ 
such that
\begin{align}\label{e31b}
 \frac{1}{\m_{2}}\G^{\frac{1}{\m_{1}}}(t-\t,X^{-1}_{\t,t}(x)-\x)&\leq \mathbf{\G}(t,x;\t,\x)\leq
 \m_{2}\G^{\m_{1}}(t-\t,X^{-1}_{\t,t}(x)-\x),\\ \label{e32b}
 \left| \p_{x_{i}}\mathbf{\G}(t,x;\t,\x)\right|&\leq
 \frac{\m_{2}}{\sqrt{t-\t}}\G^{\m_{1}}(t-\t,X^{-1}_{\t,t}(x)-\x),\\ \label{e33b}
 \left| \p_{x_{i}x_{j}}\mathbf{\G}(t,x;\t,\x)\right|
 &\leq \frac{\m_{2}}{t-\t}\G^{\m_{1}}(t-\t,X^{-1}_{\t,t}(x)-\x),
\end{align}
for every $0\leq\t<t\leq T$ and $x,\x\in \R^{d}$, where $\G^{\m}$ denotes the fundamental solution
of the heat equation $\p_{t}u_{t}(x)=\frac{\m}{2}\Delta u_{t}(x)$.
\end{theorem}
The proof of Theorem \ref{t2b} is postponed to Section \ref{mproof}.
\begin{corollary}\label{c26}
Let $u_0$ be a $\F_{0}\otimes\B$-measurable function on $\O\times\R^d$ such that $u_0(\w,\cdot)$
is continuous and non-rapidly increasing on $\R^{d}$ for a.e. $\w\in\O$. Let
$f\in\mathbf{C}^{\bar{\a}}_{0,T,\text{\rm loc}}$, for some $\bar{\a}\in(0,1)$, be non-rapidly
increasing uniformly on $[0,T]\times \R^{d}$. Then
\begin{equation}\label{e38b}
  u_{t}(x)=\int_{\R^d}\mathbf{\G}(t,x;0,\x)u_{0}(\x)d\x
  +\int_{0}^t\int_{\R^d}\mathbf{\G}(t,x;s,\x)\ff_s(\x)d\x ds
\end{equation}
is a classical solution of \eqref{SPDE} with initial value $u_0$, in the sense that $u\in
\mathbf{C}^{2}_{0,T}$ and with probability one satisfies
\begin{equation}\label{def}
 u_t(x)=u_0(x)+\int_0^t\big(\LL_su_s(x)+\ff_s(x)\big)ds+
 \int_0^t\Lc_{\s_s^k}u_s(x)dW^k_s,\qquad t\in[0,T],\ x\in\R^{d}.
\end{equation}
Such a solution is unique in the class of functions with quadratic exponential growth
{: precisely,
$u$ is the unique solution such that there exists a positive random variable $C$ such that
$\left|u_{t}(x)\right|e^{-C|x|^{2}}$ is bounded on $[\t,T]\times \R^{d}$.}
\end{corollary}

\subsection{Stochastic heat equation and Duhamel principle}\label{sec21}
To further illustrate and motivate our results, in this section we consider the prototype case of
the stochastic heat equation. We focus our attention on the Duhamel principle that is the crucial
ingredient in the parametrix method for the construction of the fundamental solution. More
generally, the Duhamel principle is a powerful tool for studying the existence and regularity
properties of parabolic PDEs. In the framework of SPDEs of the form \eqref{SPDE}, it is
still possible to have a Duhamel representation when 
the coefficients $\aa^{ij}$ are deterministic and $\Lc_{\s_t^k}\equiv 0$: this case has been
considered in \cite{Shimizu} and \cite{Mikulevicius} where the Cauchy problem for parabolic SPDEs
is studied. For the general SPDE \eqref{SPDE} however, as also noticed by other authors (see, for
instance, Sowers \cite{Sowers}, Sect.3), measurability issues arise that do not appear in the
deterministic case.

{To be more specific, let us consider the stochastic heat equation
\begin{equation} \label{ex1}
 du_t(x)=\frac{\aa^2}{2}\partial_{xx}u_t(x)dt+\left(\s \partial_x
 u_t(x)+\mathbf{g}_{t}(x)\right)dW_t.
\end{equation}
Under the coercivity condition $a^2:=\aa^{2}-\s^2>0$, the Gaussian kernel
\begin{align}\label{e10}
  p\left(t,x;\t,\x\right):=
  \frac{1}{\sqrt{2\pi a^2(t-\t) }}\exp\left( -\frac{(x+\s (W_{t}-W_{\t})-\x)^2}{2a^2(t-\t)} \right),\qquad t>\t\ge 0,\
  x,\x\in\R,
\end{align}
is well defined, and if $\s=0$ or $\mathbf{g}\equiv 0$ then the function
\begin{equation}\label{e11}
  u_{t}(x):=\int_{\R}p\left(t,x;\t,\x\right)u_{0}(\x)d\x+
  \int_{0}^{t}\int_{\R}p\left(t,x;s,\x\right)\mathbf{g}_{s}(\x)d\x\, dW_{s}
\end{equation}
is a classical solution to \eqref{ex1}, for any suitable initial value $u_{0}$. This follows
directly from the It\^o formula and the fact that the change of variable
$$X_{\t,t}(x)=x-\s(W_{t}-W_{\t})$$ transforms \eqref{ex1} into a {\it deterministic} heat
equation.}

{The difficulty in considering the case when $\s$ and $\mathbf{g}$ are {\it both} not null,
comes from the fact that the integrand $p\left(t,x;s,\x\right)\mathbf{g}_{s}(\x)$ in \eqref{e11}
becomes measurable with respect to the {\it future} $\s$-algebra $\F_{t}$ in the filtered space:
thus in general the last integral in \eqref{e11} is not well-defined in the framework of classical
It\^o-based stochastic calculus. For this reason, in the context of SPDEs, the Duhamel principle
has been used only under
rather specific assumptions.}

We observe that a naive application of the parametrix method for SPDE \eqref{SPDE} would consist
precisely of a successive application 
of the Duhamel formula \eqref{e11} with $\mathbf{g}$ and $\s=\s_{t}(x)$ that are not null and not
even constant. Hence, the lack of a general Duhamel formula seems to preclude a direct use of the
whole parametrix approach.

Incidentally formula \eqref{e10} shows that, even for SPDEs with constant coefficients, the
stochastic fundamental solution $p$ has distinctive properties compared to the Gaussian
deterministic heat kernel. In particular, the asymptotic behaviour near the pole of $p$ is
affected by the presence of the Brownian motion: this fact was studied also in \cite{Sowers} in
the more general framework of Riemannian manifolds and is coherent with the Gaussian lower and
upper bounds \eqref{e31b}.

\section{It\^o-Wentzell change of coordinates}\label{ito}
In this section we consider the  random change of coordinates \eqref{IW-SDE} and use the
It\^o-Wentzell formula to transform the SPDE \eqref{SPDE} into a PDE with random coefficients. For
simplicity, we only consider the case $\t=0$ and set $X_{t}(x)\equiv X_{0,t}(x)$. We define the
operation ``hat'' which transforms any function $u_{t}(x)$ into
\begin{equation}\label{e20}
  \hat{u}_{t}(x)=u_{t}(X_{t}(x))
\end{equation}
and recall the classical It\^o-Wentzell formula (see, for instance, Theor. 3.3.1 in
\cite{MR1070361} or Theor. 6.4 in \cite{Krylov17}).
\begin{theorem}[\bf It\^o-Wentzell] Let
$u\in \mathbf{C}^{2}_{0,T}$, $h\in \mathbf{C}^{0}_{0,T}$ and $g^k\in \mathbf{C}^{1}_{0,T}$ be such
that
\begin{equation}\label{wentzell0}
  du_{t}(x)={h_{t}(x)}dt+
  g^{k}_{t}(x)dW^{k}_{t}.
\end{equation}
Then we have
\begin{equation}\label{wentzell}
  d\hat{u}_{t}(x)=\left({\hat{h}_{t}(x)}+\frac{1}{2}\widehat{\sigma_{t}^{ik}\sigma_{t}^{jk}}(x)\widehat{\p_{ij}u_{t}}(x)
  -\widehat{\p_{i}g^{k}_{t}}(x)\hat{\s}^{ik}_{t}(x)
  \right)dt+
  \left(\hat{g}_{t}^{k}(x)-\widehat{\Lc_{\s^{k}_{t}}u_{t}}(x)\right)dW^{k}_{t}.
\end{equation}
\end{theorem}
In order to apply It\^o-Wentzell formula to our SPDE, we prove the following crucial estimate for
the gradient of $X_{t}(x)$.
\begin{proposition}\label{lemma1}
Let $$Y_{t}:=(\DD X_{t})^{-1}. 
$$ We have $\DD X_{t}$, $Y_{t}\in \mathbf{bC}^{1}_{0,T}$ and there exists a positive
random variable $\tilde{\MM}$ such that 
\begin{equation}\label{e21}
 \left|Y_t^{\ast}(x)\x\right|^{2} \geq \tilde{\MM}|\x|^2,\qquad  t\in [0,T],\ x,\x\in \R^{d},\ P\text{-a.s.}
\end{equation}
\end{proposition}
The proof of Proposition \ref{lemma1} is based on the following preliminary lemma:
%

\begin{lemma}\label{lem2}
Let $Z$ be a continuous random field defined on $[\t,T]\times\R^d$. Assume that for some $\e>0$
and $p>\left( d\vee \frac{d}{2\e}\right)$ there exists a constant $C>0$ such that
\begin{align}
 E\left[\sup_{s\in [\t,T]}\left| Z_s(x) \right|^p \right]&\leq C (1+|x|^2)^{-\e p}, \label{E2}\\
 E\left[\sup_{s\in [\t,T]}\left| \DD Z_s(x)\right|^p \right]&\leq C (1+|x|^2)^{-\e p}, \label{E3}
\end{align}
for every $x\in\R^d$. Then $Z$ has a modification in $\mathbf{bC}^{1-\frac{d}{p}}_{\t,T}$.
\end{lemma}
\proof 
By the classical Sobolev embedding theorem, for every $f\in W^{1,p}(\R^d)$, with $p>d$, we have
 $$|f(x)|+\frac{|f(x)-f(y)|}{|x-y|^{1-\frac{d}{p}}}\leq N \|f\|_{W^{1,p}(\R^d)},\qquad \text{a.e. }x,y\in\R^{d},$$
where $N$ is a constant dependent only on $p$ and $d$. Hence the statement directly follows from the
following estimate
\begin{equation}
 \sup_{t \in [\t,T]}\| Z_t \|_{W^{1,p}(\R^d)}<\infty \quad P\text{-a.e.}
\end{equation}
and the continuity of $Z$. To this end, we check that
\begin{equation}
 E\left[ \sup_{t \in [\t,T]}\| Z_t \|_{W^{1,p}(\R^d)}\right]<\infty.
\end{equation}
By \eqref{E2} and since $p>\frac{d}{2\e}$, we have
\begin{equation}
 E\left[ \sup_{t \in [\t,T]}\| Z_t \|^p_{L^{p}(\R^d)}\right]
 \leq  E\left[ \int_{\R^d}\sup_{t \in [\t,T]}| Z_t(x)|^pdx\right]
 \leq \int_{\R^d}C(1+|x|^2)^{-\e p}dx <\infty,
\end{equation}
and analogously by \eqref{E3} we have
 $$E\left[ \sup_{t \in [\t,T]}\| \DD Z_t \|^p_{L^{p}(\R^{d})}\right]
 \leq \int_{\R^d}C(1+|x|^2)^{-\e p}dx <\infty.$$
\endproof

\proof[Proof of Proposition \ref{lemma1}] Let
\begin{equation} \label{eq6}
  Z_t(x):=\DD X_t(x)-I=\int_0^t \DD\s_s^k(X_s(x))\DD X_s(x)dW^k_s.
\end{equation}
We show that the  matrix-valued random field $Z_t(x)$ satisfies estimates \eqref{E2} and
\eqref{E3} of Lemma \ref{lem2} for every $p$ such that $\left(2\vee d\vee
\frac{d}{2\e}\right)<p<\bar{p}$, with $\e$ and $\bar{p}$ as in Assumption \ref{Ass3}. Indeed, by
the well-known $L^p$-estimates for $X_{t,T}(x)$ (see \cite{MR1070361}, Chapter 4), for any
$0\le\t\le t\le T$ and $x\in \R^d$ we have
\begin{align}
 &E\left[(1+|X_{\t,t}(x)|^2)^{p} \right]\leq N_{1} (1+|x|^2)^p, \qquad p\in \R, \label{E4}\\
 &E\left[\left| \p_{}^{\b}X_{\t,t}(x)\right|^{p} \right]\leq N_{2}, \qquad p\geq 2, \ 1\le |\b|\leq 3, \label{E5}
\end{align}
where the constants $N_{1}$ and $N_{2}$ depend only on $p$ and $d$. We have
\begin{align}
E\left[\sup_{t\in [0,T]}\left| Z^{ij}_t(x) \right|^p\right]&\le C
\sum_{k=1}^{d_1}\sum_{h=1}^dE\left[\sup_{t\in [0,T]}
\left|\int_0^t\p_h\s_s^{ik}(X_s(x))\p_jX_s^h(x)dW^k_s\right|^p\right] \intertext{(by Burkolder
inequality
)}
 &\le C'_p \sum_{k=1}^{d_1}\sum_{h=1}^d
E\left[\left(\int_0^T\left(\p_h\s_s^{ik}(X_s(x))\p_jX_s^h(x)\right)^2ds\right)^{\frac{p}{2}}\right]
\intertext{(by H\"older inequality with conjugate exponents $\frac{p}{2}$ and $\frac{p}{p-2}$)}
&\le C'_p T^{\frac{p-2}{2}} \sum_{k=1}^{d_1}\sum_{h=1}^d
 \int_0^TE\left[\left|\p_h\s_s^{ik}(X_s(x))\p_jX_s^h(x)\right|^p\right]ds \intertext{(by
H\"older inequality with conjugate exponents $r$ and $q<\frac{\bar{p}}{p}$)} &\leq C'_p
T^{\frac{p-2}{2}} \sum_{k=1}^{d_1}\sum_{h=1}^d \int_0^T E\left[\left|
\p_{h}\s_s^{ik}(X_s(x))\right|^{pq}\right]^{\frac{1}{q}} E \left[\left|
\p_jX^h_s(x)\right|^{pr}\right]^{\frac{1}{r}}ds \intertext{(by Assumption \ref{Ass3} and
estimate \eqref{E5})}
&\leq C''_p T^{\frac{p-2}{2}}N_{2}^{\frac{1}{r}} \int_0^T E\left[M^{pq}(1+|X_s(x)|^2)^{-\e pq}\right]^{\frac{1}{q}}ds \\
 \intertext{(by H\"older inequality with conjugate exponents $\bar{r}$ and $\bar{q}:=\frac{\bar{p}}{pq}>1$)}
 &\leq C''_p T^{\frac{p-2}{2}}N_{2}^{\frac{1}{r}} \int_0^T E\left[M^{\bar{p}}\right]^{\frac{p}{\bar{p}}}
 E\left[(1+|X_s(x)|^2)^{-\e pq\bar{r}} \right]^{\frac{1}{q\bar{r}}}ds
 \intertext{(by estimate \eqref{E4})}
 &\leq C''_p T^{\frac{p}{2}} N^{\frac{1}{q\bar{r}}}_{1} N^{\frac{1}{r}}_{2}\|M\|_{\bar{p}}^p(1+|x|^2)^{-\e p}.
\end{align}
This proves \eqref{E2}. Estimate \eqref{E3} is obtained in a similar way from the identity
$\p_hZ^{ij}_t(x)=\p^2_{hj}X^i_t $, with $\p^2_{hj}X^i_t$ satisfying SDE \eqref{eq5b}, and
employing estimate \eqref{E5} with $|\b|=2$. Hence, by Lemma \ref{lem2}, $Z_t(x)$ has a
$\mathbf{bC}^{1-\frac{d}{p}}_{0,T}$-modification and therefore $\DD X_t(x)$ is bounded as a
function of $(t,x)\in[0,T]\times \R^{d}$, $P$-a.e. by \eqref{eq6}.

Next we prove 
that $\det \DD X_t(x)$ is bounded from above and below by a positive random variable 
for all $(t,x)$, $P$-a.s. By It\^o formula (see \cite{Krylov17}, Lemma 3.1 for more details), with
probability one we have
\begin{equation}\label{e23}
\det \DD X_t(x)=\exp\left(-\int_0^t \text{tr} D\s^k_s(X_s(x))dW^k_s
+\frac{1}{2}\int_0^t\text{tr}\left((D\s^k_s)^2\right)(X_s(x))ds \right).
\end{equation}
Since both parts of the equality are continuous w.r.t $(t,x)$, the equality holds for all
$(t,x)$ at once with probability one. Thus the assertion follows from the boundedness of the
integrals appearing in \eqref{e23}, which again can be proved as an application of Lemma
\ref{lem2}, estimate \eqref{E4} and Assumption \ref{Ass3}.

Then the matrix $Y_t(x)$ is well defined and $\det Y_t(x)$ is bounded from below by a positive random variable 
for all $(t,x)$, $P$-a.s. This fact, together with the uniform boundedness 
of the entries of $\DD X_t(x)$, implies \eqref{e21}.

It remains to prove that $\DD X_t$ and $Y_t$ have uniformly bounded spatial derivatives $P$-a.s.
Again, this is a consequence of formula \eqref{eq5b}, Lemma \ref{lem2} and the simple equality
$\p_jY_t(x)=-Y_t(x)\p_j(\DD X_t(x))Y_t(x)$.
\endproof

\begin{theorem}\label{t1}
{The function $u$ is a classical solution of SPDE \eqref{SPDE} if and only if $\hat{u}$ in
\eqref{e20} solves}
\begin{equation}\label{SPDE1}
 d\hat{u}_{t}(x)=\left(L_t\hat{u}_t(x)+f_t(x)\right)dt
\end{equation}
where $f_{t}=\hat{\ff}_{t}$ and
\begin{equation}\label{PDE}
 L_t =\frac{1}{2}
 a_t^{ij}\p_{ij}+
 b_{t}^{i} \p_i+c_t
\end{equation}
is the parabolic operator with coefficients $a^{ij},b^j,c\in\mathbf{bC}^{\a}_{0,T}$ given explicitly by 
\begin{align}\label{e25}
 a^{ij}_{t}=\left(Y_{t}\hat{\AA}_{t}Y_{t}^{\ast}\right)_{ij},\qquad
 b^{i}_{t}=
 Y^{ir}_{t}\left(\hat{\bb}^{r}_{t}-
 \widehat{\p_{j}\s^{r k}_t\s^{jk}_t}-
 \hat{\aa}^{jh}_{t}\left(Y_{t}^{\ast}(\DD^2X^{r}_{t})Y_{t}\right)_{jh}\right),\qquad
 c_{t}=\hat{\cc}_{t}.
\end{align}
Moreover, for some positive random variable $\mu$, the following coercivity condition is satisfied
\begin{equation}\label{e22}
 \langle a_{t}(x)\x,\x\rangle\ge \mu|\xi|^{2},\qquad  t\in [0,T],\ x,\x\in \R^{d},\ P\text{-a.s.}
\end{equation}
\end{theorem}
\proof By assumption, $u_{t}$ satisfies \eqref{wentzell0} with ${h_{t}}=\LL_{t}u_{t}+\ff_{t} \in
\mathbf{C}^{\a}_{0,T}$ and $g_{t}^{k}=\Lc_{\s_t^k}u_t \in \mathbf{C}^{1+\a}_{0,T}$. Thus, by the
It\^o-Wentzell formula \eqref{wentzell} we get
\begin{equation}\label{wentzell1}
 d\hat{u}_{t}=\left(\frac{1}{2}
 \hat{\AA}^{ij}_{t}\widehat{\p_{ij}u_{t}}+\left(\hat{\bb}_t^j-\widehat{\p_{i}\s^{jk}_t\s^{ik}_t}\right)
 \widehat{\p_ju}_t+\hat{\cc}_t\hat{u}_t+\hat{\ff}\right)dt.
\end{equation}
Now, we have
\begin{align}
 \p_j\hat{u}_t(x)&=
 \widehat{\p_iu}_t(x)\p_jX^i_t(x)=\left(\widehat{\DD
 u_t}(x)\DD X_t(x)\right)_j,\\
 \p_{ij}\hat{u}_t(x)&=\left(\DD X^{\ast}_t(x)\widehat{\DD^2u_t}(x)
 \DD X_t(x)\right)_{ij}+ 
 \left( \widehat{\p_hu_t}(x) \DD^2X_t^h(x) \right)_{ij},
\intertext{or equivalently}
 \widehat{\DD u_t}(x)&=\DD\hat{u}_t(x)Y_t(x),\\
 \widehat{\DD^2u_t}(x)&=Y^{\ast}_t(x)\DD^2\hat{u}_t(x)Y_t(x)-\left(Y^{\ast}_t(x)\DD^2 X^h_t(x)
 Y_t(x)\right)\widehat{\p_hu_t}(x).
\end{align}
Plugging these formulas into \eqref{wentzell1} and rearranging the indexes, we get
\eqref{SPDE1}-\eqref{PDE}-\eqref{e25}. Moreover, from expressions \eqref{e25} combined with
Assumption \ref{Ass2} and Proposition \ref{lemma1} it is straightforward to see that
$a^{ij},b^j,c\in\mathbf{bC}^{\a}_{0,T}$. Eventually, by Assumption \ref{Ass1} and estimate
\eqref{e21} of Proposition \ref{lemma1} we have
\begin{equation}
 \langle \hat{\AA}_t(x)Y^{\ast}_t\xi, Y^{\ast}_t\xi \rangle \ge \MM |Y^{\ast}_t(x)\xi|^{2}\geq \MM\tilde{\MM}|\x|^2
\end{equation}
for any $t\in [0,T]$, $x,\x\in \R^{d},$ $P$-a.s. and this proves \eqref{e22}.
\endproof

\section{Time-dependent parametrix}\label{parametrixsec}
In this section we consider the (deterministic) parabolic PDE
\begin{align}\label{PDE0}
 \Hh u_{t}(x):=L_tu_{t}(x)-\p_{t}u_{t}(x)=0
\end{align}
where
\begin{align}\label{PDE0b}
 L_tu_{t}(x)=\frac{1}{2}
 a_t^{ij}\p_{ij}u_{t}(x)+
 b_{t}^{i} \p_iu_{t}(x)+c_tu_{t}(x)
\end{align}
appears in the reduced equation \eqref{SPDE1} when $\w\in\O$ is fixed. Since the coefficients will
be assumed only measurable in the time variable, equation \eqref{PDE0} has to be understood in the
integral sense: a solution to the Cauchy problem
\begin{equation}\label{e37}
  \begin{cases}
    \Hh u_{t}(x)+f_{t}(x)=0, & x\in\R^{d},\ \text{a.e. }t\in(\t,T],\\
    u_{\t}(x)=\phi(x), & x\in\R^{d},
  \end{cases}
\end{equation}
is a function $u\in C_{\t,T}^{2}(\R^{d})$ that satisfies
\begin{equation}\label{e36}
 u_t(x)=\phi(x)+\int_{\t}^{t}(L_su_s(x)+f_s(x))ds,\qquad (t,x)\in[\t,T]\times\R^{d}.
\end{equation}

The main idea of the parametrix method is to construct the fundamental solution $\G=\G(t,x;\t,\x)$
of $\Hh $ using as a first approximation the so-called {\it parametrix}, that is the Gaussian
kernel of the heat operator obtained by freezing the coefficients of $\Hh$ at the pole $(\t,\x)$. If $Z=Z(t,x;\t,\x)$ denotes the parametrix, 
one looks for the fundamental solution of $\Hh$ in the form
\begin{equation}\label{eq0b}
 \G(t,x;\t,\x)=Z(t,x;\t,\x)+\int_{\t}^{t}\int_{\R^{d}}Z(t,x;s,y)\Phi(s,y;\t,\x)dyds.
\end{equation}
The unknown function $\Phi$ is determined by imposing $\Hh\Gamma(t,x;\t,\x)=0$: this implies that
$\Phi$ should satisfy the integral equation
\begin{equation}\label{eq0}
 \Phi(t,x;\t,\x)=\Hh Z(t,x;\t,\x)+\int_{\t}^{t}\int_{\R^{d}}\Hh Z(t,x;s,y)\Phi(s,y;\t,\x)dyds
\end{equation}
for any $x,\x\in\R^d$ and a.e. $t \in (\t,T]$. 
By recursive approximation we have
\begin{equation} \label{eq1}
 \Phi(t,x;\t,\x)=\sum_{k=1}^{+\infty}(\Hh Z)_k(t,x;\t,\x)
\end{equation}
where
\begin{align}
 (\Hh Z)_1(t,x;\t,\x)&=\Hh Z(t,x;\t,\x), \\ (\Hh Z)_{k+1}(t,x;\t,\x)&=\int_{\t}^{t}\int_{\R^{d}}\Hh Z(t,x;s,y)(\Hh Z)_k(s,y;\t,\x)dyds, \quad k\in \N.
\end{align}

To prove convergence of the series \eqref{eq1} and show that the candidate $\G$ in
\eqref{eq0b}-\eqref{eq0} is indeed a fundamental solution for $\Hh$, we need to impose some
conditions.
\begin{assumption}[\bf Coercivity]\label{ass1}
There exists a positive constant $\MMd$ such that
\begin{equation}
 \MMd^{-1}|\xi|^{2}\le \langle a_{t}(x)\x,\x\rangle\le \MMd|\xi|^{2},\qquad
 t\in[0,T],\ x,\x\in\R^{d}. 
\end{equation}
\end{assumption}
\begin{assumption}[\bf Regularity]
\label{ass2} The coefficients $a^{ij}, b^{j}, c$ are bounded functions and $a^{ij}\in
C^{\a}_{0,T}(\R^{d})$, $b^{j},c\in C^{\a}_{0,T,\text{\rm loc}}(\R^{d})$ for some $\a\in(0,1)$.
\end{assumption}
Since it is clearly not restrictive, for the subsequent analysis it is convenient to increase a
bit the value of the constant $\l$ of Assumption \ref{ass1} so that $\l>1$ and the
$L^{\infty}$-norms of the coefficients $a^{ij},b^{j},c$ are bounded by $\l$.

\begin{remark}
As opposed to the classical parametrix method, in Assumption \ref{ass2} we do not require
any regularity of the coefficients in the time variable. 
Instead, here we only require H\"older continuity in the spatial variables. The reason lies in the
fact that we are going to adopt a time-dependent definition of parametrix: namely, we do not
freeze the time variable in the definition of $Z$ (see \eqref{pareq} below) and take as parametrix
the fundamental solution of a parabolic equation with coefficients depending on $t$.
\end{remark}
\begin{remark}
Using the enhanced version of the parametrix method proposed in \cite{DeckKruse}, we can weaken
the conditions on the first- and zero-order coefficients that can be supposed to be unbounded with
sub-linear growth at infinity.
\end{remark}

\begin{definition}\label{d3} A fundamental solution $\G=\G(t,x;\t,\x)$ for equation
\eqref{PDE0} is a function defined for $0\le\t<t\le T$ and $x,\x\in\R^{d}$, such that for any
$(\t,\x)\in[0,T)\times\R^{d}$ we have:
\begin{itemize}
  \item[i)] $\G(\cdot,\cdot;\t,\x)\in C^{2}_{t_{0},T}(\R^{d})$ for any $t_{0}\in\,]\t,T[$ and 
  satisfies $\Hh \G(t,x;\t,\x)=0$ for any $x\in\R^{d}$ and a.e. $t\in(\t,T]$;
  \item[ii)] for any continuous and non-rapidly increasing function $\phi$ on $\R^{d}$
    $$\lim_{(t,x)\to(\t,\x)\atop t>\t}\int_{\R^{d}}\G(t,x;\t,y)\phi(y)dy=\phi(\x).$$
\end{itemize}
\end{definition}

Next we state the main result of this section.
\begin{theorem}[\bf Existence of the fundamental solution]\label{t2}
Under Assumptions \ref{ass1} and \ref{ass2}, there exists a fundamental solution $\G$ for equation
\eqref{PDE0}. Moreover, assume that $\phi=\phi(x)$ is continuous and non-rapidly increasing on
$\R^{d}$, and $f=f_{t}(x)$ is non-rapidly increasing uniformly on $[\t,T]\times \R^{d}$ and such
that $f\in C^{\a'}_{\t,T,\text{\rm loc}}$ for some $\a'\in(0,1)$. Then
\begin{equation}\label{e38}
  u_{t}(x)=\int_{\R^d}\G(t,x;\t,\x)\phi(\x)d\x+\int_{\t}^t\int_{\R^d}\G(t,x;s,\x)f(s,\x)d\x ds
\end{equation}
is a solution to the Cauchy problem \eqref{e37}. Such a solution is unique in the class of
functions with quadratic exponential growth (cf. Corollary \ref{c26}).
\end{theorem}

\begin{theorem}[\bf Properties of the fundamental solution]\label{t3}
Under the same assumptions of Theorem \ref{t2}, the fundamental solution $\G$ enjoys the following
properties:
\begin{itemize}
  \item[i)] $\G$ verifies the Chapman-Kolmogorov identity
\begin{align}\label{e39}
  &\G(t,x;t_{0},x_{0})=\int_{\R^{d}}\G(t,x;\t,\x)\G(\t,\x;t_{0},x_{0})d\x,\qquad t_{0}<\t<t,\ x,x_{0}\in\R^{d};
\intertext{and, if $c=c_{t}$ is independent of $x$, we have}
  \label{e40}
  &\int_{\R^{d}}\G(t,x;\t,\x)d\x=e^{\int_{\t}^{t}c_{s}ds},\qquad \t\le t\le T,\ x\in\R^{d}.
\end{align}
In particular, if $c\equiv 0$ then $\G(t,x;\t,\cdot)$ is a density;

  \item[ii)] there exist two positive constants $C_{i}=C_{i}(\MMd,\a,d,T)$, $i=1,2$,
  such that
\begin{align}\label{e31}
 \frac{1}{C_{1}}\G^{C_{2}}(t-\t,x-\x)&\leq \G(t,x;\t,\x)\leq
 C_{1}\G^{\MMd}(t-\t,x-\x),\\ \label{e32}
 \left| \p_{x_{i}}\G(t,x;\t,\x)\right|&\leq
 \frac{C_{1}}{\sqrt{t-\t}}\G^{\MMd}(t-\t,x-\x),\\ \label{e33}
 \left| \p_{x_{i}x_{j}}\G(t,x;\t,\x)\right|
 &\leq \frac{C_{1}}{t-\t}\G^{\MMd}(t-\t,x-\x),
\end{align}
for every $0\leq\t<t\leq T$ and $x,\x\in \R^{d}$, where $\G^{\MMd}$ denotes the fundamental
solution of the heat equation $\p_{t}u_{t}(x)=\frac{\MMd}{2}\Delta u_{t}(x)$.
\end{itemize}
\end{theorem}


\subsection{Preliminary Gaussian and potential estimates}
Let $A=\left(A^{ij}\right)_{1\leq i,j\leq d}$ be a constant, symmetric and positive definite
matrix. We denote by
\begin{equation} \label{Gheat}
\Gamma^{\text{\rm heat}}(A,x)=\frac{1}{\sqrt{(2\pi)^{d}\det A}}e^{-\frac{1}{2}\langle
  A^{-1}x,x\rangle},\qquad x\in\R^{d},
\end{equation}
the $d$-dimensional Gaussian kernel with covariance matrix $A$. Clearly $\Gamma^{\text{\rm heat}}$
is a smooth function and satisfies
\begin{equation}
\partial_{t}\Gamma^{\text{\rm heat}}(t A,x)=\frac{1}{2}\sum_{i,j=1}^{d}A^{ij}\p_{ij}\Gamma^{\text{\rm heat}}(t A,x),\qquad t>0,\ x\in\R^{d}.
\end{equation}

Now, we freeze the coefficients of $L_{t}$ in \eqref{PDE0b} at a fixed point $y\in\R^d$ and
consider the operator with time-dependent coefficients
\begin{equation}\label{PDEfixed}
 L_{t,y}=\frac{1}{2}
 a^{ij}_t(y)\p_{x_{i}x_{j}}
\end{equation}
acting in the $x$-variable. We denote by
\begin{equation}\label{e52}
  \G_y(t,x;\t,\x)=\G^{\text{\rm heat}}(A_{\t,t}(y), x-\x),\qquad A_{\t,t}(y):=\int_{\t}^ta_s(y) ds,
\end{equation}
the fundamental solution of $L_{t,y}-\p_{t}$. Notice that $\G_{y}$ is well defined for
$0\le\t<t\le T$ in virtue of Assumption \ref{ass1} and solves
  $$\p_{t}\G_y(t,x;\t,\x)=L_{t,y}\G_y(t,x;\t,\x)$$
for any $x,\x\in\R^{d}$ and almost every $t\in(\t,T]$. 
Finally, we define the parametrix for $\Hh$ as
\begin{equation}\label{pareq}
  Z(t,x;\t,\x)=\G_{\x}(t,x;\t,\x),\qquad 0\le \t<t\le T,\ x,\x\in\R^{d}.
\end{equation}

Hereafter $C=C\left(\cdot,\dots,\cdot\right)$ denotes a constant depending only on quantities
appearing in parentheses. In a given context the same letter will be used to denote different
constants depending on the same set of arguments. The following Gaussian estimates are standard
consequences of Assumption \ref{ass1}.
\begin{lemma}\label{p2}
We have
\begin{equation}\label{e41}
 \frac{1}{\MMd^{d}}\G^{\frac{1}{\MMd}}(t-\t,x-\x)\leq \G_{y}(t,x;\t,\x)\leq
 \MMd^{d}\G^{\MMd}(t-\t,x-\x),
\end{equation}
for any $0\leq\t<t\leq T$ and $x,\x,y \in \R^{d}$. Moreover, $\G_y(t,x;\t,\x)$ verifies the
Gaussian estimates \eqref{e32}-\eqref{e33} for some positive constant $C=C(\MMd,d)$.
\end{lemma}

\begin{proposition}\label{p3}
There exists $k_0\in \N$ such that, for every $\t \in [0,T[$ and $\x\in\R^{d}$, the series
\begin{equation}
\sum_{k=k_0}^{\infty}(\Hh Z)_k(\cdot,\cdot;\t,\x)
\end{equation}
converges in $L^{\infty}((\t,T]\times\R^{d})$. The function $\Phi$ defined by \eqref{eq1} solves
the integral equation \eqref{eq0}
and there exists a positive constant $C=C(\MMd,\a,d,T)$ such that
\begin{align}\label{eq4}
 |\Phi(t,x;\t,\x)|&\leq \frac{C}{(t-\t)^{1-\frac{\a}{2}}}\G^{\MMd}(t-\t,x-\x),\\ \label{eq4b}
 |\Phi(t,x;\t,\x)-\Phi(t,y;\t,\x)|&\leq C\frac{|x-y|^{\frac{\a}{2}}}{(t-\t)^{1-\frac{\a}{4}}}
 \left(\G^{\MMd}(t-\t,x-\x)+\G^{\MMd}(t-\t,y-\x)\right),
\end{align}
for every $x,y,\x\in\R^d$ and almost every $t \in (\t,T]$.
\end{proposition}
\proof We prove the preliminary estimate
\begin{equation}\label{eq2}
  \left|(\Hh Z)_k(t,x;\t,\x)\right|\leq \frac{M_k}{(t-\tau)^{1-\alpha k/2}}\Gamma^{\MMd}(t-\t,x-\x) \qquad 
   x,\x\in\R^d, \ \text{a.e. }t \in (\t,T],\ k\in\N,
\end{equation}
where $C=C(\MMd,\a,d,T)$ is a positive constant,
 $M_k=C^{k}\frac{\Gamma_E^k\left( \frac{\alpha}{2} \right)}{\Gamma_E\left( \frac{\alpha k}{2}
 \right)}$ and $\G_E$ is the Euler Gamma function.

For $k=1$, we have
  $$\left| \Hh Z(t,x;\t,\x)\right|= \left|(L_t -L_{t,\x})Z(t,x;\t,\x)\right|\leq I_1+I_2+I_3$$
where
\begin{align}
  I_{1}=\frac{1}{2} \big|a_t^{ij}(x)-a_t^{ij}(\x)\big|\left|\partial_{ij}Z(t,x;\t,\x)
   \right|,\quad &I_{2}=\left|b_t^{i}(x)\right|\left|\partial_{i}Z(t,x;\t,\x)\right|, \quad
  I_{3}=\left|c_t(x)Z(t,x;\t,\x)\right|.
\end{align}
{By Assumption \ref{ass2} and Lemma \ref{p2}, we have}
\begin{equation}
 I_1 \leq 
 \frac{C}{(t-\t)^{1-\frac{\a}{2}}}\left(\frac{|x-\x|}{\sqrt{t-\t}}\right)^{\a}\G^{\MMd}(t-\t,x-\x)
 \leq \frac{C}{(t-\t)^{1-\frac{\a}{2}}}\G^{\MMd+1}(t-\t,x-\x).
\end{equation}
Since the coefficients are bounded, by Lemma \ref{p2} we also have
\begin{align}
 &I_2\leq \frac{C}{\sqrt{t-\t}}\G^{\MMd}(t-\t,x-\x)\leq C
 \frac{\G^{\MMd}(t-\t,x-\x)}{(t-\tau)^{1-\frac{\a}{2}}},\qquad
 &I_3\leq C\G^{\MMd}(t-\t,x-\x), 
\end{align}
and this proves \eqref{eq2} for $k=1$. Now we assume that \eqref{eq2} holds for $k$ and prove it
for $k+1$: we have
\begin{align}
 \left| (\Hh Z)_{k+1}(t,x;\t,\x)\right| &=
 \left|\int_{\tau}^t \int_{\R^d}\Hh Z(t,x;s,y)(\Hh Z)_k(s,y;\tau,\xi)dy ds\right|\le
 \intertext{(by inductive hypothesis)}
 &\leq \int_{\tau}^t \frac{M_1}{(t-s)^{1-\alpha/2}}\frac{M_k}{(s-\tau)^{1-\alpha k/2}}
 \int_{\R^d}\Gamma^{\MMd}(t-s,x-y)\Gamma^{\MMd }(s-\t,y-\x)dy ds\le
\intertext{(by the Chapman-Kolmogorov property for $\Gamma^{\MMd }$)}
  &\leq\G^{\MMd}(t-\t,x-\x)\int_{\t}^t\frac{M_1}{(t-s)^{1-\alpha/2}}\frac{M_k}{(s-\tau)^{1-\alpha k/2}}ds
\end{align}
that yields \eqref{eq2} thanks to the well-known properties of the Gamma function. From
\eqref{eq2} we directly deduce the uniform convergence of the series and estimate \eqref{eq4}. The
proof of \eqref{eq4b} follows the same lines as in the classical case (see \cite{Friedman}, Ch.1,
Theor.7) 
and is omitted.
\endproof

We close this section by stating a generalization of a classical result about the so-called {\it
volume potential} defined as
\begin{equation}\label{e26}
  V_{f}(t,x)=\int_{t_{0}}^{t}\int_{\R^{d}}Z(t,x;\t,\x)f(\t,\x)d\x d\t,\qquad
  (t,x)\in[t_{0},T]\times \R^{d},
\end{equation}
where $Z$ denotes the parametrix. 
The proof is based on classical arguments (see \cite{Friedman}, Ch.1, Sec.3 and
\cite{IliKalOle62}) that can be applied to the time-dependent parametrix $Z$ in \eqref{pareq}
without any significant change.
\begin{lemma}\label{ll1}
Let $V_{f}$ be the volume potential in \eqref{e26} with $f\in C^{\a}_{t_{0},T,\text{\rm
loc}}(\R^{d})$, non-rapidly increasing uniformly w.r.t. $t$. Then $V_{f}\in
C^{2}_{t_{0},T}\left(\R^{d}\right)$ satisfies
\begin{align}
  \p_{x_{i}}V_{f}(t,x)&=\int_{t_{0}}^{t}\int_{\R^{d}}\p_{x_{i}}Z(t,x;\t,\x)f(\t,\x)d\x d\t,\\
  \p_{x_{i}x_{j}}V_{f}(t,x)&=\int_{t_{0}}^{t}\int_{\R^{d}}\p_{x_{i}x_{j}}Z(t,x;\t,\x)f(\t,\x)d\x d\t,\\
  \p_{t}V_{f}(t,x)&=f(t,x)+\int_{t_{0}}^{t}\int_{\R^{d}}\p_{t}Z(t,x;\t,\x)f(\t,\x)d\x d\t,
\end{align}
for any $x\in \R^{d}$ and a.e. $t\in(t_{0},T]$.
\end{lemma}

\subsection{Proof of Theorem \ref{t2}}
Let $\G=\G(t,x;\t,\x)$ be the function defined by \eqref{eq0b}-\eqref{eq1} for $0\le \t<t\le T$
and $x,\x\in\R^{d}$.
By Proposition \ref{p3}, it is clear that $\G(\cdot,\cdot;\t,\x)\in C^{0}_{\t,T}(\R^{d})$ for any
$(\t,\x)\in[0,T)\times\R^{d}$. Next, we fix $t_{0}\in(\t,t)$ and notice that by
\eqref{eq4}-\eqref{eq4b} the function $f:=\Phi(\cdot,\cdot;\t,\x)$, defined on $[t_{0},T]\times
\R^{d}$, satisfies the conditions of Lemma \ref{ll1}: hence the volume potential
  $$V_{\Phi}(t,x):=\int_{t_{0}}^{t}\int_{\R^{d}}Z(t,x;s,y)\Phi(s,y;\t,\x)dyds$$
is twice continuously differentiable in $x$ and satisfies
  $$\Hh V_{\Phi}(t,x)=\int_{t_{0}}^{t}\int_{\R^{d}}\Hh Z(t,x;s,y)\Phi(s,y;\t,\x)dyds-\Phi(t,x;\t,\x),\qquad \text{a.e. }
  t\in(t_{0},T].$$
On the other hand, we also have
  $$\Hh{\int_{\t}^{t_{0}}\int_{\R^{d}} Z(t,x;s,y)\Phi(s,y;\t,\x)dy ds=\int_{\t}^{t_{0}}\int_{\R^{d}}\Hh Z(t,x;s,y)\Phi(s,y;\t,\x)dy ds}$$
by the dominated convergence theorem. Consequently, we have
  $${\Hh \G(t,x;\t,\x)=\Hh Z(t,x;\t,\x)+\int_{\t}^{t}\int_{\R^{d}}\Hh Z(t,x;s,y)\Phi(s,y;\t,\x)dyds-\Phi(t,x;\t,\x)=0}$$
for a.e. $t \in (\t,T]$, because $\Phi$ solves equation \eqref{eq0}. 
This proves property i) of Definition \ref{d3} of fundamental solution. To prove property ii), it
suffices to notice that
\begin{equation}\label{e51}
 \int_{\R^{d}}\G(t,x;\t,y)\phi(y)dy=I_{1}(t,x,\t)+I_{2}(t,x,\t)
\end{equation}
where
\begin{align}
 \lim_{(t,x)\to(\t,\x)\atop t>\t}I_{1}(t,x,\t)&=\lim_{(t,x)\to(\t,\x)\atop t>\t}\int_{\R^{d}}Z(t,x;\t,y)\phi(y)dy=\phi(\x),\\
 \lim_{(t,x)\to(\t,\x)\atop t>\t}\left|I_{2}(t,x,\t)\right|&\le
 \lim_{(t,x)\to(\t,\x)\atop t>\t}\int_{\R^{d}}\int_{\t}^{t}
 \int_{\R^{d}}Z(t,x;s,\y)\left|\Phi(s,\y;\t,y)\phi(y)\right|d\y dsdy\le
\intertext{(by \eqref{e41}-\eqref{eq4} and since $\phi$ is non-rapidly increasing, taking $\d>0$
suitably small, with $C=C(\l,\d)$)}
 &\le\lim_{(t,x)\to(\t,\x)\atop t>\t}\int_{\R^{d}}\int_{\t}^{t}\frac{C}{(s-\t)^{1-\frac{\a}{2}}}
 \int_{\R^{d}} \G^{\MMd}(t-s,x-\y)\G^{\MMd}(s-\t,\y-y){e^{\d|y|^{2}}}d\y dsdy\\
 &\le\lim_{(t,x)\to(\t,\x)\atop t>\t}\int_{\R^{d}}\int_{\t}^{t}\frac{C}{(s-\t)^{1-\frac{\a}{2}}}
 \G^{\MMd}(t-\t,x-y){e^{\d|y|^{2}}}dsdy=0.\label{e53}
\end{align}
Finally, the standard proof of existence for the Cauchy problem (see for instance \cite{Friedman},
Ch.1, Theor.12, or \cite{IliKalOle62}) applies without modification. Uniqueness follows from the
maximum principle.

\subsection{Proof of Theorem \ref{t3}}\label{lowr}
The Chapman-Kolmogorov identity follows from uniqueness of the Cauchy problem \eqref{e37} and
representation \eqref{e38} with $f\equiv 0$ and $\phi=\G(\t,\cdot;\t_{0},\x_{0})$, for fixed
$(\t_{0},\x_{0})\in[0,\t)\times\R^{d}$. Analogously, formula \eqref{e40} follows from uniqueness
of the Cauchy problem \eqref{e37} with $f\equiv 0$ and $\phi\equiv1$.

Next we prove the Gaussian estimates for $\G$. First, we set
  $$J(t,x;\t,\x)=\int_{\t}^{t}\int_{\R^{d}}Z(t,x;s,y)\Phi(s,y;\t,\x)dyds$$
and notice that
\begin{align}
 \left|J(t,x;\t,\x)\right|&\le
 \int_{\t}^{t}\int_{\R^{d}}Z(t,x;s,y)\left|\Phi(s,y;\t,\x)\right|dyds\le
\intertext{(by \eqref{e41}, \eqref{eq4} and the Chapman-Kolmogorov property for
$\G^{\MMd}$)}\label{e42}
  &\leq\G^{\MMd}(t-\t,x-\x)\int_{\t}^t\frac{C}{(s-\tau)^{1-\frac{\a}{2}}}ds\leq C(t-\t)^{\frac{\a}{2}}\G^{\MMd}(t-\t,x-\x)
\end{align}
for some positive $C=C(\l,d,T)$. Since $\G=Z+J$, the previous estimate combined with \eqref{e41}
proves
  $$\left|\G(t,x;\t,\x)\right|\leq C_{1}\G^{\MMd}(t-\t,x-\x)$$
and in particular, the upper bound for $\G$ in \eqref{e31}. The proof of \eqref{e32}-\eqref{e33}
is similar. Notice that by the maximum principle (in the form of Lemma 5 p.43 in \cite{Friedman})
applied to $u(t,x)=\int_{\R^{d}}\G(t,x;\t,y)\phi(y)dy$, where $\phi$ is any bounded, non-negative
and continuous function, one easily infers that $\G$ is non-negative.

To prove the Gaussian lower bound we adapt a procedure due to Aronson that is essentially based on
a crucial Nash's lower bound (see \cite{FaSt}, Sect. 2). The main difference is that in our
setting we replace Nash's estimate with a bound that we directly derive from the parametrix
method. Let us first notice that, for $\l>1$, we have $\G^{\MMd}(t;x)\leq\G^{\frac{1}{\l}}(t;x)$
if $|x|\leq \r_{\l}\sqrt{t}$ where $\r_{\l}=\sqrt{\frac{\l d}{\l^2-1}\log\l}$. Thus, by
\eqref{e41} and \eqref{e42} we have
\begin{align}
 \G(t,x;\t,\x)&\geq Z(t,x;\t,\x)-|J(t,x;\t,\x)|\ge
\intertext{(if $|x-\x|\leq \r_{\l}\sqrt{t-\t}$)}
 &\geq \left(\l^{-d}-C(t-\t)^{\frac{\a}{2}}\right)\G^{\frac{1}{\l}}(t-\t;x-\x)\\ &\geq
 \frac{1}{2\l^d}\G^{\frac{1}{\l}}(t-\t;x-\x)\label{stima_loc}
\end{align}
if $0< t-\t\le T_{\l}:=(2C\lambda^d)^{-\frac{2}{\a}}\wedge T$.

For any $(t,x), (\t,\x)\in [0,T]\times\R^{d}$, we set $m$ to be the smallest natural number
greater than
\begin{equation}\label{e45}
  \max{\left\{4\r_{\l}^{-2}\frac{|x-\x|^2}{(t-\t)}, \frac{T}{T_{\l}} \right\}}.
\end{equation}
Then we set
\begin{align}\label{e0}
 t_i=\t+i\frac{t-\t}{m+1},\quad x_i=\x+i\frac{x-\x}{m+1}, \qquad i=0,\dots, m+1.
\end{align}
Denoting by $D(x;r)=\{y\in\R^{d}\mid |x-y|<r\}$ the Euclidean ball centered at $x$ with radius
$r>0$, by the Chapman-Kolmogorov equation we have
\begin{align}
 \G(t,x;\t,\x)=\int_{\R^{md}}&\G(t,x;t_{m},y_{m})\prod_{i=1}^{m-1} \G(t_{i+1},y_{i+1};t_i,y_i)\G(t_{1},y_{1};\t,\x)dy_1
 \cdots dy_{m}
\intertext{(since $\G$ is non-negative)}
 \geq \int_{\R^{md}}&\G(t,x;t_{m},y_{m})\caratt_{D(x_m;r)}(y_{m})\prod_{i=1}^{m-1}
 \G(t_{i+1},y_{i+1};t_i,y_i)\caratt_{D(x_{i};r)}(y_{i})\G(t_{1},y_{1};\t,\x)dy_1\cdots dy_{m}.\label{e47}
\end{align}
Now we have
 $$t_{i+1}-t_i=\frac{t-\t}{m+1}\leq \frac{T}{m+1}\leq T_{\l},\qquad i=0,\dots, m$$
by definition of $m$. Moreover, if $y_{i}\in D(x_{i};r)$ for $i=1,\dots, m$, by the triangular
inequality we have
\begin{align}
  \left|y_{i+1}-y_{i}\right|&\le 2r+\left|x_{i+1}-x_{i}\right|=2r+\frac{\left|x-\x\right|}{m+1}\le
\intertext{(again, by definition of $m$)}
  &\le 2r+\frac{\r_{\l}}{2}\sqrt{\frac{t-\t}{m+1}}\le
  \r_{\l}\sqrt{\frac{t-\t}{m+1}},
  \label{e46}
\end{align}
if we set
  $$r= \frac{\r_{\l}}{4}\sqrt{\frac{t-\t}{m+1}}>0.$$
For such a choice of $r$, we can use \eqref{stima_loc} repeatedly in \eqref{e47} 
and get
\begin{align}
 \G(t,x;\t,\x)&\geq (2\l^d)^{-(m+1)}
 \int_{\R^{md}}\G^{\frac{1}{\l}}\left(\frac{t-\t}{m+1},x-y_m\right)\caratt_{D(x_m;r)}(y_{m})
 \prod_{i=1}^{m-1}
 \G^{\frac{1}{\l}}\left(\frac{t-\t}{m+1},y_{i+1}-y_i\right)\times \\
 &\qquad\qquad\qquad\qquad \times\caratt_{D(x_{i};r)}(y_{i})
 \G^{\frac{1}{\l}}\left(\frac{t-\t}{m+1},y_{1}-\x\right)dy_1\cdots dy_{m}
\intertext{(by \eqref{e46} and denoting by $\w_{d}$ the volume of the unit ball in $\R^d$)}
 &\geq (2\l^d)^{-(m+1)}(\w_dr^{d})^m
 \left(\frac{\l(m+1)}{2\pi(t-\t)}\right)^{\frac{d}{2}(m+1)}\exp\left(-\frac{\l\r^2_{\l}}{2}(m+1)\right).
\end{align}
It follows that there exists a positive constant $C=C(\MMd,\a,d,T)$ such that
\begin{equation}
 \G(t,x;\t,\x)\geq \frac{1}{C(t-\t)^{\frac{d}{2}}}\exp\left(-C m\right),
\end{equation}
and this implies the required estimate.

%
%
%
%
\endproof

\subsection{Proof of Theorem \ref{t2b}}\label{mproof}
For any fixed $\tt\in[0,T)$, we consider the stochastic flow $X_{\tt,t}$ defined as in
\eqref{IW-SDE} for $t\in[\tt,T]$. Let $L^{(\tt)}_{t}$ be the operator defined as in
\eqref{PDE}-\eqref{e25} through the random change of variable $X_{\tt,t}$. By Theorem \ref{t1},
$L^{(\tt)}_{t}$ is a parabolic operator on the strip $[\tt,T]\times \R^{d}$ with random
coefficients, that satisfies Assumptions \ref{ass1} and \ref{ass2} for almost every $\w\in\O$.
Then, by Theorem \ref{t2}, $L^{(\tt)}_{t}$ admits a fundamental solution $\G^{(\tt)}(t,x;\tt,\x)$
defined for $t\in(\tt,T]$ and $x,\x\in\R^{d}$. We put
\begin{equation}\label{e54}
 \mathbf{\G}(t,x;\t,\x):=\G^{(\tt)}\left(t,X_{\tt,t}^{-1}(x);\tt,\x\right),\qquad t\in(\tt,T],\ x,\x\in\R^{d}.
\end{equation}
Combining Theorems \ref{t1} and  \ref{t2}, we infer that $\mathbf{\G}(\cdot,\cdot;\t,\x)\in
\mathbf{C}^{2}_{\t,T}(\R^{d})$ and satisfies \eqref{e50} with probability one. Moreover, let us
consider a continuous and non-rapidly increasing function $\phi$ on $\R^{d}$; proceeding as in the
proof of Theorem \ref{t2} we have
\begin{align}
  \int_{\R^{d}}\mathbf{\G}(t,x;\t,y)\phi(y)dy-\phi(\x)=I_{1}(t,x,\t)+I_{2}(t,x,\t)
\end{align}
where $I_{2}(t,x,\t)$ is defined and can be estimated as in \eqref{e53}; whereas, recalling the
definition of parametrix in \eqref{e52}\eqref{pareq}, we have
\begin{align}
 \lim_{(t,x)\to(\t,\x)\atop t>\t}I_{1}(t,x,\t)&=\lim_{(t,x)\to(\t,\x)\atop t>\t}\int_{\R^{d}}\G^{\text{\rm heat}}\left(A_{\t,t}(y),
 X_{\tt,t}^{-1}(x)-y\right)\left(\phi(y)-\phi(\x)\right)dy\\
 &=\lim_{(t,x)\to(\t,\x)\atop t>\t}\int_{\R^{d}}\G^{\text{\rm
 heat}}\left(A_{\t,t}\left(X_{\tt,t}^{-1}(x)-y\right),y\right)\left(\phi\left(X_{\tt,t}^{-1}(x)-y\right)-\phi(\x)\right)dy=0
\end{align}
by the dominated convergence theorem. This proves that $\mathbf{\G}$ is a fundamental solution for
the SPDE \eqref{SPDE}.

The Gaussian bounds \eqref{e31b} follow directly from the definition \eqref{e54} and the analogous
estimates \eqref{e31} for $\G^{(\tt)}$ in Theorem \ref{t3}. Moreover, since
  $$\p_{x_{i}}\mathbf{\G}(t,x;\t,\x)=(\nabla\G^{(\tt)})\left(t,X_{\tt,t}^{-1}(x);
  \tt,\x\right)\p_{i}X_{\tt,t}^{-1}(x),$$
the gradient estimate \eqref{e32b} follows from the analogous estimate \eqref{e32} for
$\G^{(\tt)}$ and from Proposition \ref{lemma1}. The proof of \eqref{e33b} is analogous.

{\begin{remark}
Assumption \ref{Ass3} is crucial in that it allows
us to prove 
that $\DD X_t(x)$ in \eqref{eq6} 
is finite $P$-almost surely. In turn, this guarantees that the PDE obtained through the
It\^o-Wentzell change of coordinates is uniformly parabolic. 
It would be interesting to investigate to what extent it can be relaxed and if it is possible to
have some explicit estimate.
\end{remark}}
\begin{remark}
 More general SPDEs of the form
  $$dv_{t}(x)=\left(\LL_tv_t(x)+\ff_t(x)\right)dt+
  \left(\Lc_{\s_t^k}v_t(x)+\hh^{k}_{t}(x)v_t(x)+\gg^{k}_{t}(x)\right)dW^k_t,$$
can be included in our analysis. Following the argument in \cite{MR3334279}, proof of Theor.2.7,
the idea would be to consider the transformation
 $$v_{t}(x)=\mathcal{E}_{t}(x)\left(u_{t}(x)+\int_{0}^{t}\mathcal{E}^{-1}_{s}(x)\gg^{k}_{s}(x)dW^k_s\right),$$
where
  $$\mathcal{E}_{t}(x)=\exp\left(\int_{0}^{t}\hh^{k}_{s}(x)dW^{k}_{s}-\frac{1}{2}
  \int_{0}^{t}\left|\hh^{k}_{s}(x)\right|^{2}ds\right).$$
In fact, if $u$ solves \eqref{SPDE} then by the \Ito formula we have
\begin{equation}
\begin{split}
  dv_t(x)=&\left(\mathcal{E}_{t}(x)\LL_t \left(\mathcal{E}^{-1}_{t}(x)v_t(x)\right)+\mathcal{E}_{t}(x)\ff_t(x)\right)dt\\
  &+
  \left(\mathcal{E}_{t}(x)\Lc_{\s_t^k}\left(\mathcal{E}^{-1}_{t}(x)v_t(x)\right)+\hh^{k}_{t}(x)v_t(x)+\gg^{k}_{t}(x)\right)dW^k_t.
\end{split}
\end{equation}
\end{remark}

\bibliographystyle{acm}
\bibliography{bib}

\end{document}